\documentclass{article}
\usepackage{geometry}
\usepackage{graphicx}
\usepackage{amsmath}
\usepackage{amssymb}
\usepackage{tikz}
\usetikzlibrary{shapes.geometric,arrows,positioning}
\usetikzlibrary{external}
\usepackage{multirow}
\usepackage{listings}
\usepackage{xcolor}
\definecolor{codegreen}{rgb}{0,0.6,0}
\definecolor{codegray}{rgb}{0.5,0.5,0.5}
\definecolor{codepurple}{rgb}{0.58,0,0.82}
\definecolor{backcolour}{rgb}{0.95,0.95,0.92}
\usepackage{authblk}

\usepackage[numbers,sort&compress]{natbib}
\usepackage{setspace}
\renewcommand{\vec}{\boldsymbol}
\newcommand{\x}{\vec{x}}
\newcommand{\y}{\vec{y}}
\newcommand{\xivec}{\vec{\xi}}
\newcommand{\R}{\mathbb{R}}

\newcommand{\z}{\vec{y}}
\newcommand{\U}{\mathcal{U}}

\lstdefinestyle{python}{
    backgroundcolor=\color{backcolour},
    commentstyle=\color{codegreen},
    keywordstyle=\color{magenta},
    numberstyle=\tiny\color{codegray},
    stringstyle=\color{codepurple},
    basicstyle=\ttfamily\footnotesize,
    breakatwhitespace=false,
    breaklines=true,
    captionpos=b,
    keepspaces=true,
    numbers=left,
    numbersep=5pt,
    showspaces=false,
    showstringspaces=false,
    showtabs=false,
    tabsize=2
}

\begin{document}

\title{ROmodel: Modeling robust optimization problems in Pyomo
}
\author[1]{Johannes Wiebe}
\author[1]{Ruth Misener}
\affil[1]{Department of Computing, Imperial College London, London, UK}
\date{}

\maketitle

\begin{abstract}
    This paper introduces ROmodel, an open source Python package extending the
    modeling capabilities of the algebraic
    modeling language Pyomo to robust optimization problems.
    ROmodel helps practitioners transition from deterministic to robust
    optimization through modeling objects which allow formulating robust models
    in close analogy to their mathematical formulation.
    ROmodel contains a library of commonly used uncertainty sets which can be
    generated using their matrix representations, but it also allows users to
    define custom uncertainty sets using Pyomo constraints. ROmodel
    supports adjustable variables via linear decision rules. The
    resulting models can be solved using ROmodels solvers which implement both
    the robust reformulation and cutting plane approach.
    ROmodel is a platform to implement and compare custom uncertainty
    sets and reformulations. We demonstrate ROmodel's capabilities by applying
    it to six case studies. We implement custom
    uncertainty sets based on (warped) Gaussian processes to show how ROmodel
    can integrate data-driven models with optimization.
\end{abstract}

\section{Introduction}
\label{intro}
Robust optimization is a common way of managing optimization under uncertainty
in process systems engineering: applications range from production scheduling
to flexible chemical process design
\cite{Janak2005,Li2008,Zhang2015,Ning2017,Shang2018,Grossmann2016}.
New developments in robust optimization include
distributionally and adjustable robust optimization to reduce solution conservatism
\cite{Grossmann2016},
data-driven robust optimization to model uncertainty sets
based on available data \cite{Bertsimas2018}, and approximate robust optimization
to make non-linear problems more tractable \cite{Houska2013}.
This large number of alternative approaches and the required domain knowledge
can discourage practitioners from making the transition from deterministic
optimization to optimization under uncertainty. Furthermore, the lack of a
platform that allows the easy implementation and application of new algorithms
means that it is difficult for researchers to compare different
approaches \cite{Goerigk2016}.

This paper introduces ROmodel, a
Python package extending the popular, Python-based algebraic modeling language Pyomo
\cite{hart2011pyomo,Hart2017} to facilitate modeling of robust optimization
problems and implementation of robust optimization algorithms. ROmodel aims to
(i) make robust optimization more accessible to practitioners and facilitate
moving from deterministic to robust optimization, (ii)
demonstrate how robust optimization problems can be modeled in close analogy to
their mathematical formulation, and (iii) provide an open-source platform for
researchers to implement and compare new robust reformulations and uncertainty
sets. To this end, ROmodel introduces intuitive modeling
which closely resembles the optimization problem's
underlying mathematical formulation and will feel familiar to Pyomo users.
ROmodel uses the richness of Pyomo's solver interfaces and
methods for model transformations and Python's data processing capabilities.
It supports both automatic reformulation and cutting plane algorithms.
Uncertainty sets can be chosen from a
library of common geometries, or custom defined using Pyomo constraints.
ROmodel also support adjustable robust optimization through linear decision
rules. Because ROmodel and Pyomo are open-source, ROmodel can be extended to
incorporate additional uncertainty set geometries and reformulations. As an
example, we implemented uncertainty sets based on (warped) Gaussian processes
for black-box constrained problems.

There are other software packages for \emph{solving} robust optimization
problems, e.g. ROME\cite{Goh2011}, RSOME\cite{Chen2020},
ROC++\cite{Vayanos2020}, and PyROS\cite{pyros}.
All four approaches are designed as robust optimization \emph{solvers} and,
except for PyROS which is based on Pyomo, do not rely
on a general purpose algebraic modeling language.
In contrast, ROmodel focuses on \emph{modeling} robust optimization problems.
By building on Pyomo, ROmodel simplifies the transition from
deterministic to robust optimization, since both the deterministic and robust
model can be
implemented in the same environment. ROmodel also allows
access to a much larger number of deterministic solvers than these other solvers.
In contrast to AIMMS, which has some capabilities for modeling robust
optimization problems \cite{aimms}, ROModel is open source and allows
possibilities for extension.
JumPeR \cite{Dunning2016}, which extends Julia's modeling language JumP to
robust optimization problems, is the most similar to ROModel. In contrast to
JumPeR, ROModel is more tightly tied to its respective algebraic modeling
language, e.g. developing new convex uncertainty sets in ROmodel only requires
adding Pyomo constraints while JumPeR would require adding Julia functionality.
ROmodel also automatically recognizes uncertainy set geometries and applies the
corresponding transformations without the user having to specify which
geometry they are using.

A further advantage of ROmodel is that it is based on
Python, which is popular in data analytics and machine learning. ROmodel therefore allows
data-based techniques to be integrated seemlesly with robust optimization
methods.
As an example, we implement (warped) Gaussian process-based uncertainty sets in
ROmodel\cite{Wiebe2020}.
 These sets seemlesly integrate Gaussian processes
trained in the Python library GPy \cite{gpy} into robust optimization problems.
ROmodel is open source and available on Github \cite{romodel}.

The rest of this paper is structured as follows. Section~\ref{sec:modeling}
introduces ROmodel's new modeling objects and shows how they can be used to
model robust optimization problems. Section~\ref{sec:solvers} introduces
ROmodel's three
solvers: a reformulation based solver, a cutting plane
solver, and a nominal solver for obtaining nominal solutions of robust
problems. Section~\ref{sec:gp} discusses how ROmodel can be extended and
demonstrates our
implementation of Gaussian process-based uncertainty sets for black-box
constrained problems. Section~\ref{sec:results} introduces six case studies
and presents results.

\section{Modeling}
\label{sec:modeling}

Consider the following generic deterministic optimization problem
\begin{subequations}
    \begin{align}
        \min\limits_{\x\in \mathcal{X} , \y}\quad & f(\x, \y, \bar{\xivec})  \\[4pt]
        \text{s.t} \quad & g(\x, \y, \bar{\xivec}) \leq 0
    \end{align}\label{eq:nominal}
\end{subequations}
where $\x \in \R^n$ and $\y \in R^m$ are decision variables and
$\bar{\xivec} \in \R^k$ is a vector of (nominal) parameters.
If the parameter vector $\bar{\xivec}$ is not known exactly, we can construct
the following robust version of Problem~\ref{eq:nominal}:
\begin{subequations}
    \begin{align}
        \min\limits_{\x\in \mathcal{X} , \y(\xivec)}\max\limits_{\xivec \in \U(\x)}\quad & f(\x, \y(\xivec), \xivec)  \\[4pt]
        \text{s.t} \quad & g(\x, \y(\xivec), \xivec) \leq 0 && \forall \xivec \in \U(\x)
    \end{align}\label{eq:ro}
\end{subequations}
Here $\x \in \R^n$ are "here and now" variables, determined before the uncertainty is revealed,
while $\y(\xivec)$ are adjustable variables, determined after the uncertainty is revealed.
The uncertain parameter vector $\xivec$ is bounded
by the uncertainty set $\U(\x)$, which may depend on $\x$ and which contains the nominal values $\bar{\xivec}$ from
Problem~\ref{eq:nominal}.
Optimization Problem~\ref{eq:ro}, a generic robust optimization problem with uncertainty in
the objective and constraints, is the basis for ROmodel.
Note that we are limiting ourselves here to one
uncertain constraint for simplicity of notation only, ROmodel can handle
multiple robust constraints. Also note that there can be an arbitrary number of deterministic
constraints which define the set $\mathcal{X}$.

ROmodel introduces three new modeling objects to represent robust optimization problems like Problem~\ref{eq:ro} within Pyomo:
\begin{enumerate}
    \item \lstinline{UncParam}: A class similar to Pyomo's \lstinline{Param}
        and \lstinline{Var} class used to model uncertain parameters $\xivec$.
        One can supply a \lstinline{nominal}
        argument, which defines the vector of nominal values $\bar{\xivec}$
        used to replace the uncertain parameters when solving the deterministic
        (nominal) problem (Eq.~\ref{eq:nominal}).
    \item \lstinline{UncSet}:
        Each \lstinline{UncParam} object has an \lstinline{UncSet}
        object associated with it.
        The \lstinline{UncSet} class, based on Pyomo's \lstinline{Block} class,
        models the uncertainty sets $\U$.  Uncertainty sets $\U$ can be
        defined by (i) constructing generic sets by adding Pyomo
        constraints to the \lstinline{UncSet} object, or (ii) through a library of
        common uncertainty set geometries, using their matrix representation as
        an input.
    \item \lstinline{AdjustableVar}: A class similar to Pyomo's \lstinline{Var}
        class used to model adjustable variables which can be
        determined after some of the uncertainty has been resolved, i.e.
        $\y(\xivec)$.
\end{enumerate}

These three new modeling objects are sufficient for modeling quite generic robust
optimization problems. They are set up to allow modeling problems in an
intuitive way which is closely related to their mathematical formulation
(Problem~\ref{eq:ro}). We discuss each modeling object and how it is used to
construct robust optimization problems in the subsequent sections.

\subsection{Uncertain Parameters}
Indexed uncertain parameters are constructed in analogy to Pyomo's
\lstinline{Var} type for variables:
\begin{lstlisting}[language=Python]
    m.c = UncParam(range(3), nominal=[0.1, 0.2, 0.3], uncset=m.U)
\end{lstlisting}
The \lstinline{nominal} argument specifies a list of nominal values
$\bar{\xivec}$ for the
uncertain paramters $\xivec$. The \lstinline{uncset} argument specifies the
uncertainty set to use for these parameters. The two approaches for constructing
the uncertainty set \lstinline{m.U} are dicussed in the next two chapters.

\subsection{Generic uncertainty sets}
\label{sec:generic}
Generic uncertainty sets are constructed with the \lstinline{UncSet} class.
This class inherits from Pyomo's \lstinline{Block} class.
Users can construct generic uncertainty sets by adding Pyomo constraints to an
\lstinline{UncSet} object in the same way as they would add constraints to a
\lstinline{Block} object in Pyomo.
The following example shows how a polyhedral set can be modeled using the
\lstinline{UncSet} class:
\begin{lstlisting}[language=Python]
    from romodel import UncSet, UncParam
    # Define uncertainty set
    m.U = UncSet()
    # Define uncertain parameter
    m.w = UncParam(range(2), uncset=m.U, nominal=[0.5, 0.5])
    # Add constraints to uncertainty set
    m.U.cons1 = Constraint(expr=m.w[0] + m.w[1] <= 1)
    m.U.cons2 = Constraint(expr=m.w[0] - m.w[1] <= 1)
    ...
\end{lstlisting}
The \lstinline{uncset} and \lstinline{nominal} arguments define the uncertainty
set and a vector of nominal values for the uncertain parameter \lstinline{m.w}.
ROmodel's strategy of modeling uncertainty sets using Pyomo's \lstinline{Constraint}
modeling object is analogous to typical robust optimization formulations.
ROmodel can therefore be used to define any
uncertainty set which can be expressed in Pyomo.
However, not every set that
can be modeled with Pyomo constraints can necessarily also be solved.
Section~\ref{sec:solvers} discusses which types of uncertainty sets can be
solved using which solver in.
Users can also define multiple uncertainty sets and replace one by another:
\begin{lstlisting}[language=Python]
    # Define second uncertainty set
    m.U2 = ro.UncSet()
    # Swap uncertainty sets and solve
    m.coef.uncset = m.U2
    solver = pe.SolverFactory('romodel.cuts')
    solver.solve(m)
\end{lstlisting}

\subsection{Library uncertainty sets}
\label{sec:library}
For commonly used, standard uncertainty sets, the generic
approach (Section~\ref{sec:generic}) is
unnecessarily complicated. Therefore, ROmodel implements custom classes which can
define an uncertainty set using its matrix representation.
ROmodel implements polyhedral and ellipsoidal sets.
The user can define polyhedral sets of the form $\mathcal{U} = \left\{w \; | \; Pw \leq b\right\}$
by passing the matrix $P$ and the right hand side $b$ to the
\lstinline{PolyhedralSet} class:
\begin{lstlisting}[language=Python]
    from romodel.uncset import PolyhedralSet
    # Define polyhedral set
    m.U = PolyhedralSet(mat=[[ 1,  1],
                             [ 1, -1],
                             [-1,  1],
                             [-1, -1]],
                        rhs=[1, 1, 1, 1])
\end{lstlisting}
ROmodel creates ellipsoidal sets of the form $(w-\mu)\Sigma^{-1}(w-\mu) \leq 1$
using the \lstinline{EllipsoidalSet} class, the covariance matrix
$\Sigma$ and the mean vector $\mu$:
\begin{lstlisting}[language=Python]
    from romodel.uncset import EllipsoidalSet
    # Define ellipsoidal set
    m.U = EllipsoidalSet(cov=[[1, 0, 0],
                              [0, 1, 0],
                              [0, 0, 1]],
                         mean=[0.5, 0.3, 0.1])
\end{lstlisting}
In Section~\ref{sec:gp} we discuss how additional library sets can be added to ROmodel using data-driven uncertainty sets based on (warped) Gaussian
processes as an example.

\subsection{Constructing uncertain constraints}
After defining uncertain parameters and an uncertainty set, the user can
construct uncertain constraints implicitly by using the uncertain parameters in
a Pyomo constraint.
Consider the following deterministic Pyomo constraint:
\begin{lstlisting}[language=Python]
    # deterministic
    m.x = Var(range(3))
    c = [0.1, 0.2, 0.3]
    m.cons = Constraint(expr=sum(c[i]*m.x[i] for i in m.x) <= 0)
\end{lstlisting}
If the coefficients $c$ are uncertain, we can model the robust constraint
$c^T x \leq 0, \; \forall c \in \mathcal{U}$ as:
\begin{lstlisting}[language=Python]
    # robust
    m.x = Var(range(3))
    m.c = UncParam(range(3), nominal=[0.1, 0.2, 0.3], uncset=m.U)
    m.cons = Constraint(expr=sum(m.c[i]*m.x[i] for i in m.x) <= 0)
\end{lstlisting}
The uncertain parameter \lstinline{m.c} can be used in Pyomo constraints in the
same way as Pyomo \lstinline{Var} or \lstinline{Param} objects. ROmodel's
solvers automatically recognize constraints and objectives containing
containing uncertain parameters.

\subsection{Adjustable variables}
ROmodel also has capabilities for modeling adjustable variables. Adjustable
variables $\y \left(\xivec\right)$ are variables which can be determined after some (or all) of the
uncertainty has been revealed. Defining an adjustable variable is analogous to defining a
regular variable in Pyomo, with an additional \lstinline{uncparam} argument
specifying a list of uncertain parameters which the adjustable variable depends
on:
\begin{lstlisting}[language=Python]
    # Define uncertain parameters
    m.w = UncParam(range(3), nominal=[1, 2, 3])
    # Define adjustable variable which depends on uncertain parameter
    m.y = AdjustableVar(range(3), uncparams=[m.w], bounds=(0, 1))
\end{lstlisting}
The uncertain parameters can also be set individually for each element of the
adjustable variables index using the \lstinline{set_uncparams} function:
\begin{lstlisting}[language=Python]
    # Set uncertain parameters for individual indicies
    m.y[0].set_uncparams([m.w[0]])
    m.y[1].set_uncparams([m.w[0], m.w[1]])
\end{lstlisting}
ROmodel only implements linear decision rules for solving adjustable robust
optimization problems. If a model contains
adjustable variables in a constraint or objective, ROmodel automatically
replaces it by a linear decision rule based on the specified uncertain
parameters.

\section{Solvers}
\label{sec:solvers}
ROmodel has three solvers: A robust reformulation based solver,
a cutting plane based solver, and a nominal solver.

\subsection{Reformulation}

\begin{figure}[htb!]
    \centering
    \includegraphics{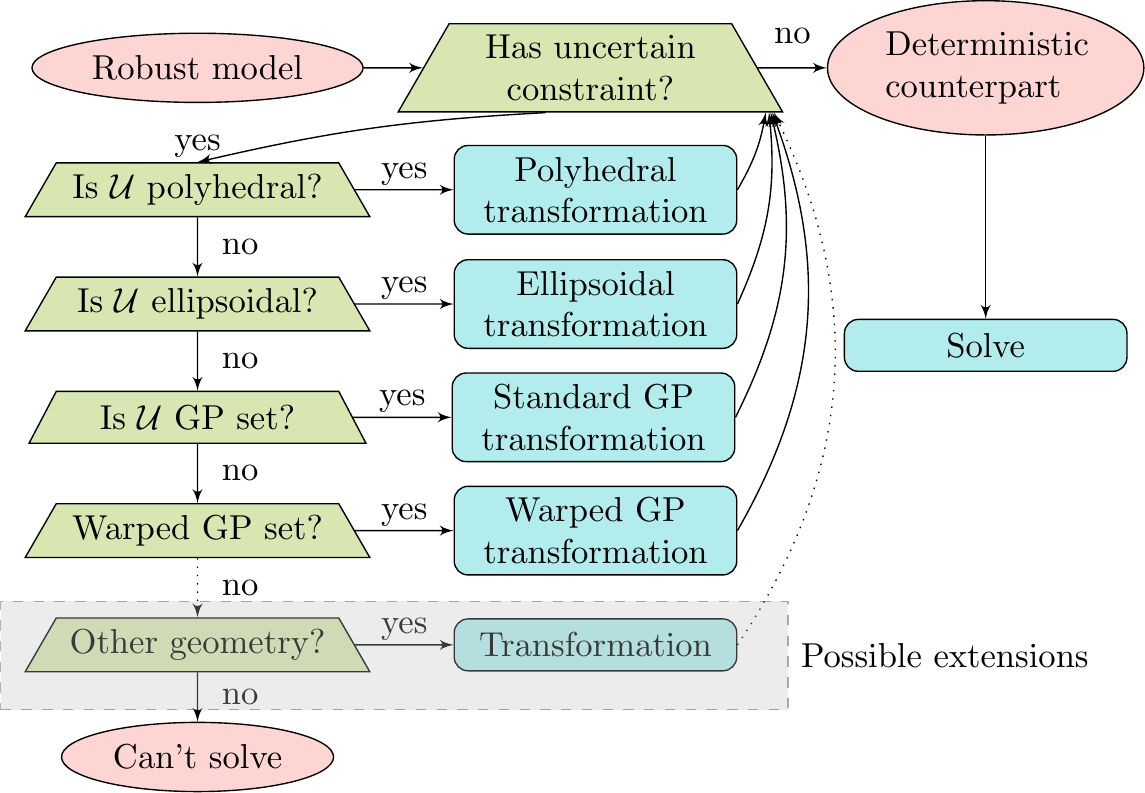}
    \caption{Schematic of reformulation solver. For each uncertain constraint,
        the reformulation solver goes
        through all known transformations and identifies the geometry of the
        constraint/uncertainty set. It then applies the corresponding model
        transformation. If all uncertain constraints can be reformulated, the
        resulting deterministic counterpart is solved using one of Pyomo's solver
        interfaces. If one or more constraints cannot be reformulated, the
        problem cannot be solved.}
    \label{fig:reform}
\end{figure}

The reformulation-based solver, illustrated in Fig.~\ref{fig:reform}, implements
standard duality based techniques for
reformulating robust optimization problems into deterministic
counterparts \cite{Bertsimas2004}.
First, it detects every
constraint containing uncertain parameters. Second, it checks the structure of
each uncertain constraint and the corresponding uncertainty set to determine if
a known reformulation is applicable. Finally, it applies a model
transformation, generating the deterministic counterpart of each robust
constraint. The deterministic counterpart is then solved using an
appropriate solver available in Pyomo. The structure of the optimization
problem and the uncertainty set geometry determine which solvers are
applicable. If no applicable transformation can be identified for one or more
constraints, the problem cannot be solved and ROmodel will raise an error.

ROmodel implements standard reformulations for ellipsoidal and
polyhedral uncertainty sets and linear uncertain constraints \cite{Bental1999,Bertsimas2004}. It also
implements reformulations for black-box constrained problems \cite{Wiebe2020}. These are
discussed in more detail in Section~\ref{sec:gp}, which also dicussed how
ROmodel can be extended to include further reformulations.

\subsection{Cutting planes}

\begin{figure}[htb!]
    \centering
    \includegraphics{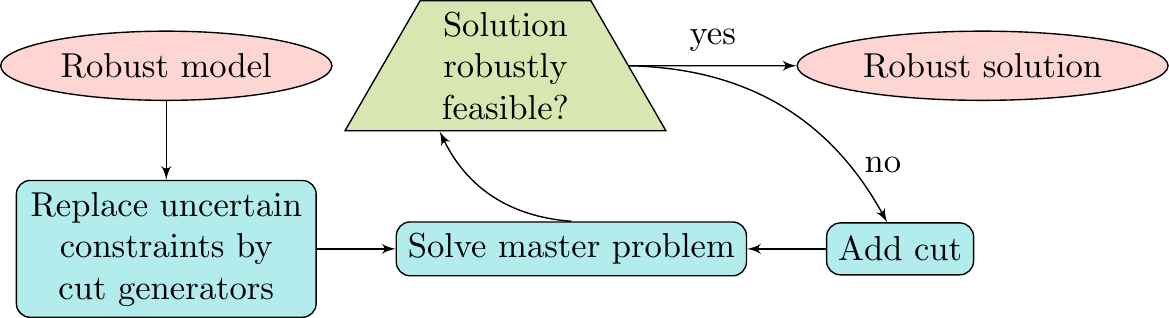}
    \caption{Schematic of cutting plane solver. The cutting plane solver
    iteratively solves the master problem and adds cutting planes until a
    robustly feasible solution is found.}
    \label{fig:cuts}
\end{figure}

The cutting plane solver, outlined in Fig.~\ref{fig:cuts}, implements an
iterative strategy for solving robust optimization problems \cite{Mutapcic2009}.
It replaces each uncertain constraint and objective by a
\lstinline{CutGenerator} object which initially just contains the nominal
constraint. The solver then iteratively solves the master problem and
generates cuts to cut off solutions which are not robustly feasible.
A solution $x^*$ is considered to be robustly feasible when for each uncertain
constraint $g(\x^*, \y(\xivec), \xivec) \leq 0$, the objective value of the
separation problem is smaller than some tolerance $\epsilon$:
    \begin{align}
        \max\limits_{\xivec \in \U} \; & g(\x^*, \y(\xivec), \xivec) \leq
        \epsilon
    \end{align}
ROmodel's cutting plane solver can generally be applied to any convex
uncertainty set. The Pyomo solvers for solving the master and separation
problems can be set individually using:
\begin{lstlisting}[language=Python]
    solver = SolverFactory('romodel.cuts')
    # Master solver
    solver.options['solver'] = 'gurobi'
    # Sub solver
    solver.options['subsolver'] = 'ipopt'
\end{lstlisting}
The solvers need to be appropriate for the corresponding problem, i.e., the
above choice would be sensible if the master problem is a mixed-integer linear
problem and the uncertain constraints and uncertainty set are continuous and
convex. If the solvers are not appropriate, Pyomo will raise an error.

\subsection{Nominal}
    ROmodel also includes a nominal solver. This solver replaces all
    occurrences of the uncertain parameters by their nominal values and solves
    the resulting deterministic problem:

    \begin{lstlisting}[language=Python]
        # Obtaining the nominal solution
        solver = SolverFactory('romodel.nominal')
        solver.solve(m)
    \end{lstlisting}
    The nominal solver allows users to combine their implementations of the
    nominal and robust problem. An implementation of the robust model can
    be used to obtain
    the solution of the nominal problem.

\section{Extending ROmodel for black-box constrained problems}
\label{sec:gp}
ROmodel can be extended to incorporate additional reformulations and
uncertainty set geometries.
This section outlines how ROmodel can be extended using (warped) Gaussian
process-based uncertainty sets for black-box constrained problems as an
example.
This example showcases the ease with which ROmodel can integrate Python's
machine learning and data analytics capabilities with Pyomo's mathematical
optimization modeling.

Wiebe et al. \cite{Wiebe2020} propose a robust optimization-based approximation
of a class of chance-constraints containing uncertain black-box functions
$g(\cdot)$:
\begin{equation}
    P\left(
        \sum\limits_i g(\z_i)x_i \leq b\right) \geq 1 - \alpha
    \label{eq:ccgp}
\end{equation}
The approach models the black-box function and associated uncertainty using a
(warped) Gaussian process as a stochastic model. The standard Gaussian process
is well known and commonly used as a surrogate model \cite{Bhosekar2018}. Warped Gaussian processes
are a more flexible variant of standard Gaussian process in which observations are mapped
into a latent space using a non-linear, often neural net-style warping function
\cite{Snelson2003}.
If a standard Gaussian process models $g(\cdot)$, the chance
constraint Eq.~\ref{eq:ccgp} can be reformulated exactly. For the warped
Gaussian process, Wiebe et al.~\cite{Wiebe2020} propose an
approximation based on Wolfe duality.

In order to make these approaches available in ROmodel, we need to (i)
implement two library uncertainty sets, \lstinline{GPSet} and
\lstinline{WarpedGPSet},  which collect the relevant data, and
(ii) implement two corresponding model transformations which perform the reformulations
for standard and warped Gaussian process-based sets.
The implementation is based on the Python module ROGP \cite{rogp}, which is
includes Gaussian process models trained in the Python library GPY
\cite{gpy} in Pyomo models.

\subsection{Implementing new library sets}
Implementing a new library set mainly requires a new
Python class collecting the necessary data.
For the standard and warped Gaussian process set, this data consist of
three arguments:
\begin{lstlisting}[language=Python]
    from romodel.uncset import GPSet, WarpedGPSet
    # Define variables
    m.y = Var([0, 1])
    # Define uncertainty set based on standard GP
    m.uncset_standard = GPSet(gp_standard, m.z, 0.95)
    # Define uncertainty set based on warped GP
    m.uncset_warped = WarpedGPSet(gp_warped, m.z, 0.95)
\end{lstlisting}
    The first argument \lstinline{gp_standard/warped} is a (warped) Gaussian
    process object trained in GPy. The second is an indexed Pyomo variable on
    which the GP depends, i.e. $\y$ in Eq.~\ref{eq:ccgp}.
    The third parameter specifies the confidence level $1-\alpha$ with which
    the true parameter is contained in the uncertainty set. I.e., in this case
    the confidence that the true parameter vector is an element of the
    uncertainty set is at least 95\%.

The new sets can be used in the same way as other library sets, e.g.:
\begin{lstlisting}[language=Python]
    # Define variables
    m.x = Var([0, 1])
    # Define uncertain parameters
    m.w = UncParam([0, 1], uncset=m.uncset_warped)
    # Define constraint
    m.cons = Constraint(expr=m.x[0]*m.w[0] + m.x[1]*m.w[1] <= 1)
\end{lstlisting}
Constraints which use this type of uncertainty set need to be linear in the
uncertain parameter. Note that the indices of \lstinline{m.z} and
\lstinline{m.w} need to be identical in the formulation above.
If the black-box function depends on
more than one variable, the Gaussian process-based sets can alternatively
be specified using a dictionary:
\begin{lstlisting}[language=Python]
    # Define variables
    m.y = Var([0, 1], ['a', 'b'])
    # Define uncertainty set
    y_dict = {0: [m.y[0, 'a'], m.y[0, 'b']],
              1: [m.y[1, 'a'], m.y[1, 'b']]},
    m.uncset_warped = WarpedGPSet(gp, y_dict, 0.95}
\end{lstlisting}
    The dictionary indicates that the uncertain parameter \lstinline{m.w[0]}
    depends on the variables
    \lstinline{m.z[0, 'a']} and \lstinline{m.z[0, 'b']} through the black-box
    function $g(\cdot)$, modeled in GPy by the Gaussian process \lstinline{gp}.

    Note that ROmodel's cutting plane solver is not applicable to the Gaussian
    process-based sets because the sets are decision dependent. Attempting to solve
    a problem with one of theses sets therefore results in an error.
    When implementing new library sets which \emph{can} be solved
    using cutting planes, an additional Python function
    \lstinline{generate_cons_from_lib}, which generates
    Pyomo constraints for the uncertainty set based on the data collected by
    the library set, is required. For an example, see
    \lstinline{romodel/uncset/ellipsoidal.py} on the ROmodel Github
    \cite{romodel}.

\subsection{Implementing new reformulations}
For ROmodel to be able to solve models containing the two new Gaussian
process-based sets, we need to
implement the corresponding reformulations.
Adding new reformulations to ROmodel generally requires
two Python functions: (i) a function \lstinline{_check_applicability} which detects
whether a constraint and uncertainty set have the required structure,
and (ii) a function \lstinline{_reformulate} which generates the
robust counterpart. The former function is only required if the reformulation
is supposed to work with generically constructed uncertainty sets as described
in Section~\ref{sec:generic}. For library sets like the Gaussian process-based
sets, only the latter function is required. This function takes data describing
the constraints and uncertainty set as an input and returns a Pyomo block
containing the deterministic counterpart.
For a full example see the implemented reformulations in
\lstinline{romodel/reformulate/} on the ROmodel Github \cite{romodel}.

\section{Results}
\label{sec:results}
We use ROmodel to model and solve six case studies:
\begin{enumerate}
    \item A portfolio optimization problem with uncertain returns \cite{Bertsimas2004},
    \item A knapsack problem with uncertain item weights,
    \item A pooling problem instance \cite{adhya} with uncertain product demands,
    \item A capacitated facility location problem as an example for adjustable
        robust optimization, where the decision which facilities to build
        has to be made under demand uncertainty, while the decision from which
        facility to supply individual customers can be made once the
        uncertainty is resolved,
    \item A production planning in which the price at which products can be
        sold depends on the amount produced through an uncertain black-box
        function modelled by a (warped) Gaussian process \cite{Wiebe2020},
    \item And a drill scheduling problem in which the equipment used to drill a
        well degrades at a rate which depends other drill parameters through a
        black-box function \cite{Wiebe2020}.
\end{enumerate}
All examples except for the drill scheduling example are included
with ROmodel and can be used as follows:

    \begin{lstlisting}[language=Python]
        import romodel.example as ex
        portfolio = ex.Portfolio()
        knapsack = ex.Knapsack()
        pooling = ex.Pooling()
        facility = ex.Facility()
        planning = ex.ProductionPlanning(alpha=0.95, warped=True)
    \end{lstlisting}
The implementation of the drill scheduling example is separately available on
Github \cite{drilling}.
We solve the portfolio, knapsack, pooling, and facility location problems with
both the reformulation and cutting
plane solver for ellipsoidal and polyhedral uncertainty sets and using both the
library approach to generating uncertainty sets as well as the generic, Pyomo
constraint-based approach. We solve the production planning and drill
scheduling problems using the reformulation solver with uncertainty sets based
on both standard and warped Gaussian processes. We solve 30 instances with
different uncertainty set sizes for each case study.

\begin{table}[htb!]
    \centering
\begin{tabular}{l l c c c}
                               &             & Reformulation    & Cuts & Overall \\ \hline
    \multirow{2}{*}{Knapsack}  & Polyhedral  & \phantom{601}54  & \phantom{40}272 & \phantom{753}85      \\
                               \cline{2-5}
                               & Ellipsoidal & \phantom{601}50  & \phantom{40}183 & \phantom{753}91      \\
    \hline
    \multirow{2}{*}{Pooling}   & Polyhedral  & \phantom{601}74  & \phantom{40}329 & \phantom{75}173     \\
                               \cline{2-5}
                               & Ellipsoidal & \phantom{60}638  & \phantom{40}331 & \phantom{75}349    \\
    \hline
    \multirow{2}{*}{Portfolio} & Polyhedral  & \phantom{601}50  & \phantom{40}276 & \phantom{75}126      \\
                               \cline{2-5}
                               & Ellipsoidal & \phantom{601}49  & \phantom{4}1659 & \phantom{75}129      \\
    \hline
    \multirow{2}{*}{Facility}  & Polyhedral  & \phantom{60}261  & 13353           & \phantom{7}5588\\
                               \cline{2-5}
                               & Ellipsoidal & \phantom{601}--  & 31275           & 31275\\
    \hline
    \multirow{2}{*}{Planning}  & Standard    & \phantom{6}2776  & \phantom{400}NA & \phantom{7}2776\\
                               \cline{2-5}
                               & Warped      & \phantom{6}8536  & \phantom{400}NA & \phantom{7}8536\\
    \hline
    \multirow{2}{*}{Drilling}  & Standard    & 13646            & \phantom{400}NA & 13646 \\
                               \cline{2-5}
                               & Warped      & 75325            & \phantom{400}NA & 75325\\
    \hline
                        Overall&             & \phantom{601}74  & \phantom{40}330 & \phantom{75}271     \\
    \hline
\end{tabular}
\caption{Median time in milliseconds taken to solve the six example problems
    with different uncertainty set geometries. Times are shown for the reformulation
    and cutting plane solvers and both combined (Overall). The cutting plane solver is not applicable (NA) for
Gaussian process-based sets and the reformulation solver cannot solve the
facility problem with ellipsoidal sets to optimality (--).}
\label{tab:results}
\end{table}
Table~\ref{tab:results} shows the median time in milliseconds taken to solve
each problem for a given uncertainty set geometry and solver. The reformulation
solver generally outperforms the cutting plane solver with median times of 74ms
and 330ms respectively. An exception is the
the non-linear, non-convex pooling problem with an ellipsoidal set. For this
instance, the cutting plane solver achieves significantly better results, which
is in line with previous work on robust pooling problems \cite{Wiebe2019}. Similarly, for the facility
location problem with an ellipsoidal set, the reformulation approach does not
solve the problem to optimality within a 10 minute time frame, while the cutting
plane solver does.
For the production planning and drills scheduling examples only the
reformulation solver can be applied. The Wolfe duality-based reformulation for
warped Gaussian processes generally takes longer to solve than the
chance-constraint reformulation for standard Gaussian processes.
Note that most of this time is the time taken by the
subsolvers. The transformations which ROmodel performs are generally very
quick: the median transformation time across all instances is $3.2$
milliseconds, while the maximum transformation time is $1.7$ seconds.

\begin{figure}[htb!]
    \centering
    \includegraphics{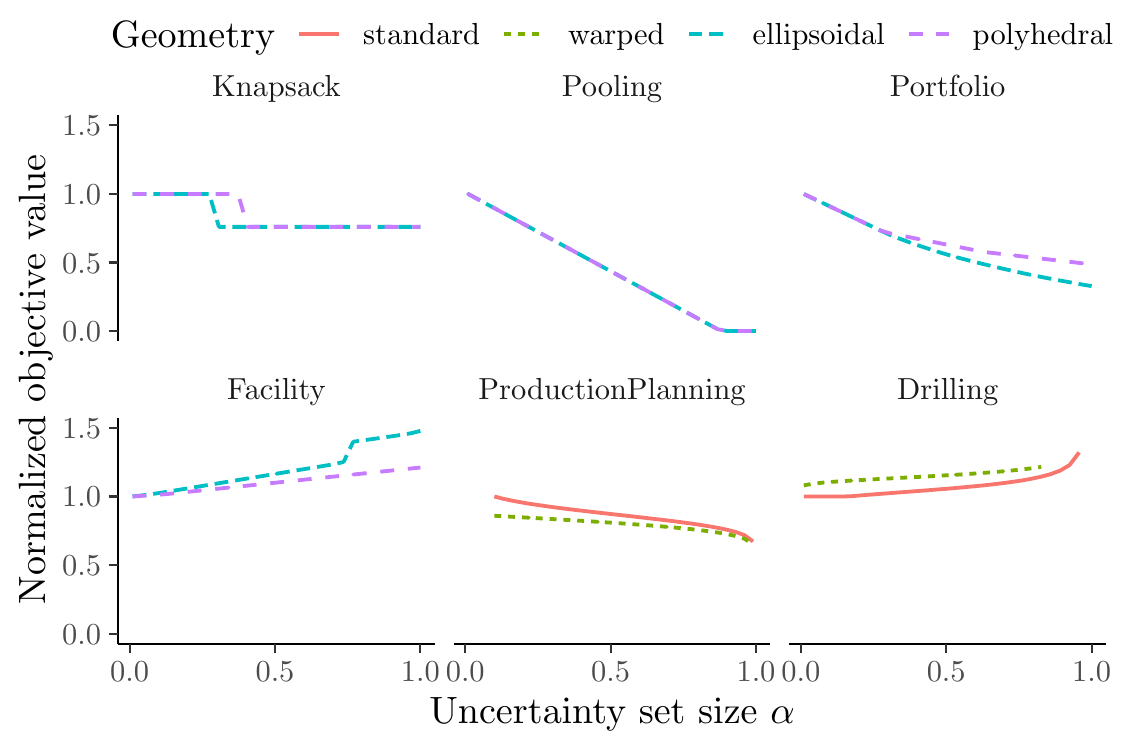}
    \caption{The figure shows the normalized objective value as a function of
    the uncertainty set size for knapsack, portfolio and pooling case studies
    and an ellipsoidal and polyhedral uncertainty set.}
    \label{fig:results}
\end{figure}
Fig.~\ref{fig:results} shows the objective value (normalized with the objective
of the nominal solver) as a function of uncertainty set size for each of the
six examples. Note that the facility location and drill scheduling problems are
minimization problems, while the rest are maximization problems. By
construction, the ellipsoidal sets tested always fully contain the polyhedral
sets for a given $\alpha$. Correspondingly, they are always more conservative,
i.e., for a given $\alpha$ the objective value achieved using the ellipsoidal
set is larger for minimization and smaller for maximization problems than the
value achieved using the polyhedral set. For the Gaussian process-based sets,
the standard approach is always less conservative than the warped approach.
However, the limited ability of the standard Gaussian process to model
non-Gaussian noise may mean that the actual probability of constraint violation
is larger than the intended confidence would suggest.
For a more detailed comparison of these two approaches see \cite{Wiebe2020}.

\section{Conclusion}

ROmodel formulates robust versions of common optimization
problems. The modeling environment it provides makes (adjustable) robust optimization
methods more readily available to practitioners and makes trying different
solution approaches and uncertainty sets very easy.
ROmodel is open source and available free of charge and could play a vital role
as a platform for prototyping
novel robust optimization algorithms and comparing them to existing approaches.

\section{Acknowledgements}
This work was funded by the Engineering \& Physical Sciences Research Council (EPSRC) Center for Doctoral Training in High Performance Embedded and Distributed Systems (EP/L016796/1), an EPSRC/Schlumberger CASE studentship to J.W. (EP/R511961/1, voucher 17000145), and an EPSRC Research Fellowship to R.M. (EP/P016871/1).

\bibliography{lit.bib}   

\end{document}